\documentclass[runningheads,a4paper]{llncs}

\usepackage{amssymb,amsmath}
\usepackage{graphicx}
\usepackage{algorithm}
\usepackage{algpseudocode}
\usepackage{color}
\usepackage{subfigure}
\usepackage[toc,page]{appendix}
\usepackage{setspace}
\usepackage{url}

\begin{document}

\mainmatter  

\title{Extended Formulations in Mixed-integer Convex Programming}

\titlerunning{Extended Formulations in Mixed-integer Convex Programming}

%
%
\author{Miles Lubin\inst{1} \and Emre Yamangil\inst{2} \and Russell Bent\inst{2} \and Juan Pablo Vielma\inst{1}}
\authorrunning{Lubin, Yamangil, Bent, Vielma}

\institute{Massachusetts Institute of Technology, Cambridge, MA, USA\\
\and
Los Alamos National Laboratory, Los Alamos, NM, USA
}

%
%

\toctitle{Lecture Notes in Computer Science}
\tocauthor{Authors' Instructions}
\maketitle

\begin{abstract}
We present a unifying framework for generating extended formulations for the polyhedral outer approximations used in algorithms for mixed-integer convex programming (MICP). Extended formulations lead to fewer iterations of outer approximation algorithms and generally faster solution times. First, we observe that 
all MICP instances from the MINLPLIB2 benchmark library are conic representable with standard symmetric and nonsymmetric cones. Conic reformulations are shown to be effective extended formulations themselves because they encode separability structure. For mixed-integer conic-representable problems, we provide the first outer approximation algorithm with finite-time convergence guarantees, opening a path for the use of conic solvers for continuous relaxations. We then connect the popular modeling framework of disciplined convex programming (DCP) to the existence of extended formulations independent of conic representability. We present evidence that our approach can yield significant gains in practice, with the solution of a number of open instances from the MINLPLIB2 benchmark library.
\end{abstract}

\section{Introduction}

Mixed-integer convex programming (MICP) is the class of problems where one seeks to minimize a convex objective function subject to convex constraints and integrality restrictions on the variables.
MICP is less general than mixed-integer nonlinear programming (MINLP), where the objective and constraints may be nonconvex, but unlike the latter, one can often develop finite-time algorithms to find a global solution. These finite-time algorithms depend on convex nonlinear programming (NLP) solvers to solve continuous subproblems. MICP, also called \textit{convex MINLP}, has broad applications and is supported in various forms by both academic solvers like Bonmin~\cite{Bonmin} and SCIP~\cite{scip} and commercial solvers like KNITRO~\cite{knitro}; see Bonami et al.~\cite{BonamiReview,MahajanReview} for a review.

The most straightforward approach for MICP is NLP-based branch and bound, an extension of the branch and bound algorithm for mixed-integer linear programming (MILP) where a convex NLP relaxation is solved at each node of the branch and bound tree~\cite{Gupta85}. However, driven by the availability of effective solvers for linear programming (LP) and MILP, it was observed in the early 1990s by Leyffer and others~\cite{SvenThesis} that it is often more effective to avoid solving NLP relaxations when possible in favor of solving polyhedral relaxations using MILP. Polyhedral relaxations form the basis of the majority of the existing solvers recently reviewed and benchmarked by Bonami et al.~\cite{BonamiReview}.

While traditional MICP approaches construct polyhedral approximations in the original space of variables, a number of authors have considered introducing auxiliary variables and forming a polyhedral approximation in a higher dimensional space~\cite{Baron,Hijazi,VielmaExtendedFormulations,MustafaThesis}. Such constructions are called \textit{extended formulations} or \textit{lifted formulations}, the motivation for which is the fact that the projection of these polyhedra onto the original space can provide a higher quality approximation than one built from scratch in the original space. Tawarmalani and Sahinidis~\cite{Baron} propose, in the context of nonconvex MINLP, extended formulations for compositions of functions. For MICP, Hijazi et al.~\cite{Hijazi} demonstrate the effectiveness of extended formulations in the special case where all nonlinear functions can be written as a sum of \textit{univariate} convex functions. Their method obtains promising speed-ups over Bonmin on the instances which exhibit this structure. Hijazi et al. generated these extended formulations by hand, and no subsequent work has proposed techniques for off-the-shelf MICP solvers to detect and exploit separability. Building on these results, Vielma et al.~\cite{VielmaExtendedFormulations} propose extended formulations for second-order cones. These extended formulations improved solution times for mixed-integer second-order cone programming (MISOCP) over state of the art commercial solvers CPLEX and Gurobi quite significantly; both solvers adopted the technique within a few months after its publication. 

A major contribution of this work is to propose a new, unifying framework for generating extended formulations for the polyhedral outer approximations used in MICP algorithms. This framework generalizes the work of Hijazi et al.~\cite{Hijazi} which was specialized for separable problems to include all MICPs whose objective and constraints can be expressed in closed algebraic form. We begin in Section~\ref{sec:conicextended} by considering conic representability. While many MICP instances are representable by using MISOCP, reformulation to MISOCP has not been widely adopted, and MICP is still considered a significantly more general form. We demonstrate that with the introduction of the nonsymmetric exponential and power cones, surprisingly, all convex instances in the MINLPLIB2 benchmark library~\cite{MINLPLIB} are representable by a combination of these nonsymmetric cones and the second-order cone. We discuss how the conic-form representation of a problem is itself a strong extended formulation. Hence, the guideline to ``just solve the conic form problem'' is surprisingly effective.

We note that conic-form problems have modeling strength beyond that of smooth MICP, in particular for handling of nonsmooth perspective functions useful in disjunctive convex programming~\cite{perspective}. With the recent development of conic solvers supporting nonsymmetric cones~\cite{akle,SCS}, it may be advantageous to use these solvers over derivative-based NLP solvers, in which case the standard convergence theory for outer approximation algorithms no longer applies. In Section~\ref{sec:conicoa}, we present the first finite-time outer approximation algorithm applicable to mixed-integer conic programming with any closed, convex cones (symmetric and nonsymmetric), so long as conic duality holds in strong form. This algorithm extends the work of Drewes and Ulbrich~\cite{MISOCPOA} for MISOCP with a much simpler and more general proof.

In Section~\ref{sec:dcp}, we generalize the idea of extended formulations through conic representability by considering the modeling framework of \textit{disciplined convex programming} (DCP)~\cite{DCP}, a popular modeling paradigm for convex optimization which has so far received little notice in the MICP realm. In DCP, convex expressions are specified in an algebraic form such that convexity can be verified by simple composition rules. We establish a 1-1 connection between these rules for verifying convexity and the existence of extended formulations. Hence, all MICPs expressed in mixed-integer disciplined convex programming (MI\textbf{D}CP) form have natural extended formulations regardless of conic representability. This view has connections with techniques for nonconvex MINLP, where it is already common practice to construct extended outer approximations based on the algebraic representation of the problem~\cite{MahajanReview}.

In our computational experiments, we translate MICP problems from the MINLPLIB2 benchmark library into MIDCP form and demonstrate significant gains from the use of extended formulations, including the solution of a number of open instances. Our open-source solver, \textit{Pajarito}, is the first solver specialized for MIDCP and is accessible through existing DCP modeling languages.

\section{Extended formulations and conic representability}\label{sec:conicextended}

We state a generic mixed-integer convex programming problem as
\begin{align}\label{eq:micvx}
\operatorname{minimize}_{x,y}\,\, & f(x,y)  \notag \\
\text{subject to } & g_j(x,y) \le 0 \quad \forall j \in J,\tag{MICONV}\\
&L \le x \le U,\quad x \in \mathbb{Z}^n, y \in \mathbb{R}^{p}_+, \notag
\end{align}
\noindent
where the set $J$ indexes the nonlinear constraints, the functions $f, g_j : \mathbb{R}^n \times \mathbb{R}^p \to \mathbb{R} \cup \{\infty \}$ are convex, and the vectors $L$ and $U$ are finite bounds on $x$. Without loss of generality, when convenient, we may assume that the objective function $f$ is linear (via epigraph reformulation~\cite{BonamiReview}).

Vielma et al.~\cite{VielmaExtendedFormulations} discuss the motivation for extended formulations in MICP: many successful MICP algorithms use polyhedral outer approximations of nonlinear constraints, and polyhedral outer approximations in a higher dimensional space can often be much stronger than approximations in the original space. Hijazi et al.~\cite{Hijazi} give an example of an approximation of an $\ell_2$ ball in $\mathbb{R}^n$ which requires $2^n$ tangent hyperplanes in the original space to prove that the intersection of the ball with the integer lattice is in fact empty. By exploiting the summation structure in the definition of the $\ell_2$ ball, \cite{Hijazi} demonstrate that an extended formulation requires only $2n$ hyperplanes to prove an empty intersection. More generally, \cite{Hijazi,Baron} propose to reformulate constraints with separable structure $\sum_{k=1}^q g_k(x_k) \le 0,$
where $g_k : \mathbb{R} \to \mathbb{R}$ are \textit{univariate} convex functions by introducing auxiliary variables $t_k$ and imposing the constraints
\begin{gather}\label{eq:hijazi}
\sum\nolimits_{k=1}^q t_k \le 0, g_k(x_k) \le t_k \forall\, k.
\end{gather}

A consistent theme in this paper is the \textit{representation} of the convex functions $f$ and $g_j,\; \forall j\in J$ in~\eqref{eq:micvx}. Current MICP solvers require continuous differentiability of the nonlinear functions and access to black-box oracles for querying the values and derivatives of each at any given point $(x,y)$. The difficulty in the reformulation~\eqref{eq:hijazi} is that the standard representation~\eqref{eq:micvx} does not encode the necessary information, since separability is an algebraic property which is not detectable given only oracles to evaluate function values and derivatives. As such, we are not aware of any off-the-shelf MICP solver which exploits this special-case structure, despite the promising experimental results of~\cite{Hijazi}.

In this section, we  consider the equally  general, yet potentially more useful, representation of~\eqref{eq:micvx} as a mixed-integer conic programming problem:
\begin{align}\label{eq:miconic}
\min_{x,z} \quad& c^Tz\notag\\
\text{s.t. } & A_xx + A_zz = b\tag{MICONE}\\
& L \le x \le U, x \in \mathbb{Z}^n, z \in \mathcal{K},\notag
\end{align}
where $\mathcal{K}\subseteq \mathbb{R}^k$ is a closed convex cone. Without loss of generality, we assume integer variables are not restricted to cones, since we may introduce corresponding continuous variables by equality constraints. The representation of~\eqref{eq:micvx} as~\eqref{eq:miconic} is equally as general in the sense that given a convex function $f$, we can define a closed convex cone $\mathcal{K}_f = \operatorname{cl}\{ (x,y,\gamma,t) : \gamma f(x/\gamma,y/\gamma) \le t, \gamma > 0 \}$ where $\operatorname{cl} S$ is defined as the closure of a set $S$. Using this, we can reformulate \eqref{eq:micvx} to the equivalent optimization problem
\begin{align}\label{eq:miconic_all}
\min \quad& t_f \notag\\
\text{s.t. } & t_j + s_j = 0 \quad \forall j \in J,\\
&\gamma_f = 1, x = x_f, y = y_f,\notag\\
&\gamma_j = 1, x = x_j, y = y_j, \forall j \in J,\notag\\
& L \le x \le U, x \in \mathbb{Z}^n, y \in \mathbb{R}^p_+,\notag\\
&(x_f,y_f,\gamma_f,t_f) \in \mathcal{K}_f\notag,\\
&(x_j,y_j,\gamma_j,t_j) \in \mathcal{K}_{g_j}, s_j \in \mathbb{R}_+ \quad\forall j \in J.\notag
\end{align}

The problem~\eqref{eq:miconic_all} is in the form of~\eqref{eq:miconic} with $\mathcal{K} = \mathbb{R}^{n+|J|}_+ \times \mathcal{K}_f \times \mathcal{K}_{g_1} \times \cdots \times \mathcal{K}_{g_{|J|}}$. Such a tautological reformulation is not particularly useful, however. What \textit{is} useful is a reformulation of~\eqref{eq:micvx} into~\eqref{eq:miconic} where the cone $\mathcal{K}$ is a product $\mathcal{K}_1 \times \mathcal{K}_2
 \times \cdots \times \mathcal{K}_r$, where each $\mathcal{K}_i$ is one of a small number of recognized cones, such as the positive orthant $\mathbb{R}^n_+$, the second-order cone $SOC_n = \{ (t,x) \in \mathbb{R}^n : ||x|| \le t \}$,
 the exponential cone, $EXP = \operatorname{cl}\{ (x,y,z) \in \mathbb{R}^3: y\exp(x/y) \le z, y > 0 \}$, and the power cone (given $0 < \alpha < 1$), $POW_\alpha = \{ (x,y,z) \in \mathbb{R}^3 : |z| \le x^\alpha y^{1-\alpha}, x \ge 0, y \ge 0\}$. 
 
 The question of which functions can be represented by second-order cones has been well studied~\cite{SOCPApplications,lectures}. More recently, a number of authors have considered nonsymmetric cones, in particular the exponential cone, which can be used to model logarithms, entropy, logistic regression, and geometric programming~\cite{akle}, and the power cone, which can be used to model $p$-norms and powers~\cite{powercone}.
 
 The folklore within the conic optimization community is that almost all convex constraints which arise in practice are representable by using these cones\footnote{\url{http://erlingdandersen.blogspot.com/2010/11/which-cones-are-needed-to-represent.html}}, in addition to the positive semidefinite cone which we do not consider here. To substantiate this claim, we classified the 333 MICP instances in MINLPLIB2 according to their conic representability and found that \textit{all} of the instances are conic representable; see Table~\ref{tab:conic}.
 
 \begin{table}[t]
 \centering
 \begin{tabular}{c|c|c|c|c|c}
   SOC only & EXP only & SOC and EXP & POW only & Not representable & Total  \\
  \hline
  217 & 107 & 7 & 2 & 0 & 333
 \end{tabular}
 \vspace{0.3cm}
 \caption{A categorization of the 333 MICP instances in the MINLPLIB2 library according to conic representability. Over two thirds are pure MISOCP problems and nearly one third is representable by using the exponential (EXP) cone alone. All instances are representable by using standard cones. }
 \label{tab:conic}
 \end{table}
 
 While solvers for SOC-constrained problems (SOCPs) are mature and commercially supported, the development of effective and reliable algorithms for handling exponential cones and power cones is an emerging, active research area~\cite{akle,SCS}. Nevertheless, we claim that the conic view of~\eqref{eq:micvx} is useful \textit{even lacking} reliable solvers for continuous conic relaxations.
 
 As a motivating example, we consider the trimloss~\cite{Harjunkoski} (\texttt{tls}) instances from MINLPLIB2, a convex formulation of the cutting stock problem. These instances are notable as some of the few unsolved instances in the benchmark library and also because they exhibit a separability structure more general than what can be handled by Hizaji et al.~\cite{Hijazi}.

The trimloss instances have constraints of the form
\begin{equation}\label{eq:tls}
\sum\nolimits_{k=1}^q -\sqrt{x_k y_k} \le c^Tz + b,
\end{equation}
where $x,y,z$ are arbitrary variables unrelated to the previous notation in this section. Harjunkoski et al.~\cite{Harjunkoski} obtain these constraints from a clever reformulation of nonconvex bilinear terms. The function $-\sqrt{xy}$ is the negative of the geometric mean of $x$ and $y$. It is convex for nonnegative $x$ and $y$ and its epigraph $E = \{(t,x,y) : -\sqrt{xy} \le t, x \ge 0, y \ge 0 \}$ is representable as an affine transformation of the three-dimensional second-order cone $SOC_3$~\cite{lectures}. A conic formulation for~\eqref{eq:tls} is constructed by introducing an auxiliary variable for each term in the sum plus a slack variable, resulting in the following constraints:
\begin{equation}\label{eq:tls1}
\sum\nolimits_{k=1}^q t_k +s = c^Tz + b,  (t_k,x_k,y_k) \in E\,\, \forall k, \text{ and } s \in \mathbb{R}_+.
\end{equation}
Equation~\eqref{eq:tls1} provides an \textit{extended} formulation of the constraint~\eqref{eq:tls}, that is, an equivalent formulation using additional variables.

If we take the MINLPLIB2 library to be representative, then conic structure using standard cones exists in the overwhelming majority of MICP problems in practice. This observation calls for considering~\eqref{eq:miconic} as a standard form of MICP, one which is perhaps more useful for computation than~\eqref{eq:micvx} precisely because it is an extended formulation which encodes separability structure in a natural and general way. There is a large body of work and computational infrastructure for automatically generating the conic-form representation given an algebraic representation, a discussion we defer to Section~\ref{sec:dcp}.

The benefits of reformulation from~\eqref{eq:micvx} to~\eqref{eq:miconic} are quite tangible in practice. By direct reformulation from MICP to MISOCP, we were able to solve to global optimality the trimloss \texttt{tls5} and \texttt{tls6} instances from MINLPLIB2 by using Gurobi 6.0\footnote{Solutions reported to Stefan Vigerske, October 5, 2015}. These instances from this public benchmark library had been unsolved since 2001, perhaps indicating that the value of conic formulations is not widely known.

\section{An outer-approximation algorithm for mixed-integer conic programming}\label{sec:conicoa}

Although the conic representation~\eqref{eq:miconic} does not preclude the use of derivative-based solvers for continuous relaxations, derivative-based nonlinear solvers are typically not appropriate for conic problems because the nonlinear constraints which define the standard cones have points of nondifferentiability~\cite{Noam}. Sometimes the nondifferentiability is an artifact of the conic reformulation (e.g., of smooth functions $x^2$ and $\exp(x)$), but in a number of important cases the nondifferentiability is intrinsic to the model and provides additional modeling power. Nonsmooth perspective functions, for example, which are used in disjunctive convex programming, have been particularly challenging for derivative-based MICP solvers and have motivated smooth approximations~\cite{perspective}. On the other hand, conic form can handle these nonsmooth functions in a natural way, so long as there is a solver capable of solving the continuous conic relaxations.

There is a growing body of work as well as some (so far) experimental solvers supporting mixed second-order and exponential cone problems~\cite{akle,SCS}, which opens the door for considering conic solvers in place of derivative-based solvers. To the best of our knowledge, however, no outer-approximation algorithm or finite-time convergence theory has been proposed for general mixed-integer conic programming problems of the form~\eqref{eq:miconic}.

In this section, we present the first such algorithm for~\eqref{eq:miconic} with arbitrary closed, convex cones. This algorithm generalizes the work of Drewes and Ulbrich~\cite{MISOCPOA} for MISOCP with a much simpler proof based on conic duality. In stating this algorithm, we hope to motivate further development of conic solvers for cones beyond the second-order and positive semidefinite cones. 

We begin with the definition of dual cones.

\begin{definition}
Given a cone $\mathcal{K}$, we define $\mathcal{K}^* := \{ \beta \in \mathbb{R}^k : \beta^Tz \ge 0 \,\, \forall z \in \mathcal{K}\}$ as the dual cone of $\mathcal{K}$.
\end{definition}

Dual cones provide an equivalent outer description of any closed, convex cone, as the following lemma states. We refer readers to~\cite{lectures} for the proof.

\begin{lemma}
Let $\mathcal{K}$ be a closed, convex cone. Then $z \in \mathcal{K}$ iff $z^T\beta \ge 0\,\, \forall \beta \in \mathcal{K}^*$.
\end{lemma}
Based on the above lemma, we will consider an outer approximation of~\eqref{eq:miconic}:
\begin{align}\label{eq:mioa}
\min_{x,z} \quad& c^Tz\notag\\
\text{s.t. } & A_xx + A_zz = b\tag{MIOA(T)}\\
& L \le x \le U, x \in \mathbb{Z}^n,\notag\\
&\beta^Tz \ge 0\,\, \forall \beta \in T.\notag
\end{align}

Note that if $T = \mathcal{K}^*$, \ref{eq:mioa} is an equivalent semi-infinite representation of~\eqref{eq:miconic}. If $T \subset \mathcal{K}^*$ and $|T| < \infty$ then~\ref{eq:mioa} is an MILP outer approximation of~\eqref{eq:miconic} whose objective value is a lower bound on the optimal value of~\eqref{eq:miconic}.

The outer approximation (OA) algorithm is based on iteratively building up $T$ until convergence in a finite number of steps to the optimal solution. First, we define the continuous subproblem for fixed integer value $\hat x$ which plays a key role in the OA algorithm:
\begin{align}\label{eq:conic_cont}
v_{\hat x} = \min_{z} \quad&  c^Tz\notag\\
\text{s.t. } & A_zz = b - A_x \hat x\tag{$CP(\hat x)$},\\
&z \in \mathcal{K}\notag.
\end{align}
The dual of~\eqref{eq:conic_cont} is
\begin{align}\label{eq:conic_cont_dual}
\max_{\beta,\lambda} \quad& \lambda^T(b-A_x \hat x)\notag\\
\text{s.t. } & \beta = c - A_z^T\lambda\\
&\beta \in \mathcal{K}^*\notag.
\end{align}

The following lemmas demonstrate, essentially, that the dual solutions to~\eqref{eq:conic_cont} provide the only elements of $\mathcal{K}^*$ that we need to consider.

\begin{lemma}\label{lem:conic}
Given $\hat x$, assume~\ref{eq:conic_cont} is feasible and strong duality holds at the optimal primal-dual solution $(z_{\hat x},\beta_{\hat x},\lambda_{\hat x})$. Then for any $z$ with $A_z z = b - A_x \hat x$ and $\beta_{\hat x}^Tz \ge 0$, we have $c^Tz \ge v_{\hat x}$.
\begin{proof}
\begin{equation}
\beta_{\hat x}^Tz = (c - A_z^T\lambda_{\hat x})^Tz = c^Tz - \lambda_{\hat x}^T(b - A_x\hat x) = c^Tz - v_{\hat x} \ge 0.
\end{equation}
\end{proof}
\end{lemma}

\begin{lemma}\label{lem:ray}
Given $\hat x$, assume~\ref{eq:conic_cont} is infeasible and~\eqref{eq:conic_cont_dual} is unbounded, such that we have a ray $(\beta_{\hat x},\lambda_{\hat x})$ satisfying $\beta_{\hat x} \in \mathcal{K}^*$, $\beta_{\hat x} = -A_z^T\lambda_{\hat x}$, and $\lambda_{\hat x}^T(b - A_x\hat x) > 0$. Then for any $z$ satisfying $A_z z = b - A_x \hat x$ we have $\beta_{\hat x}^Tz < 0$.
\begin{proof}
\begin{equation}
\beta_{\hat x}^Tz = -\lambda_{\hat x}^TA_zz = -\lambda_{\hat x}^T(b - A_x\hat x) < 0.
\end{equation}

\end{proof}

\end{lemma}

\begin{algorithm}[ht]\small
\caption{The conic outer approximation (OA) algorithm}\label{alg:oa}
\begin{algorithmic}
\State \textbf{Initialize:} $z_U \leftarrow \infty, z_L \leftarrow -\infty$, $T \leftarrow \emptyset$. Fix convergence tolerance $\epsilon$.
\While{$z_U - z_L \ge \epsilon$}
\State Solve \ref{eq:mioa}.
\If{\ref{eq:mioa} is infeasible}
\State \eqref{eq:miconic} is infeasible, so terminate.
\EndIf
\State Let $(\hat x,\hat z)$ be the optimal solution of \ref{eq:mioa} with objective value $w_T$.
\State Update lower bound $z_L \leftarrow w_T$.
\State Solve~\ref{eq:conic_cont}.
\If{\ref{eq:conic_cont} is feasible}
\State Let $(z_{\hat x},\beta_{\hat x},\lambda_{\hat x})$ be an optimal primal-dual solution with objective value $v_{\hat x}$.
\If{$v_{\hat x} < z_U$}
\State $z_U \leftarrow v_{\hat x}$
\State Record $(\hat x, z_{\hat x})$ as the best known solution.
\EndIf
\ElsIf{\ref{eq:conic_cont} is infeasible}
\State Let $(\beta_{\hat x},\lambda_{\hat x})$ be a ray of~\eqref{eq:conic_cont_dual}.
\EndIf
\State $T \leftarrow T \cup \{\beta_{\hat x}\}$
\EndWhile
\end{algorithmic}
\end{algorithm}

Finite termination of the algorithm is guaranteed because integer solutions $\hat x$ cannot repeat, and only a finite number of integer solutions is possible.

This algorithm is arguably incomplete because the assumptions of Lemmas \ref{lem:conic} and \ref{lem:ray} need not always hold. The assumption of strong duality at the solution is analogous to the constraint qualification assumption of the NLP OA algorithm~\cite{SvenOA}. Drewes and Ulbrich~\cite{MISOCPOA} describe a procedure in the case of MISOCP to ensure finite termination if this assumption does not hold. The assumption that a ray of the dual exists if the primal problem is infeasible is also not always true in the conic case, though \cite{lectures} provide a characterization of when this can occur. These cases will receive full treatment in future work.

A notable difference between the conic OA algorithm and the standard NLP OA algorithm is that there is no need to solve a second subproblem in the case of infeasibility, although some specialized NLP solvers may also obviate this need~\cite{Filmint}. In contrast, Drewes and Ulbrich~\cite{MISOCPOA} propose a second subproblem in the case of MISOCP even when dual rays would suffice.

Finally, the algorithm is presented in terms of a single cone $\mathcal{K}$ for simplicity. When $\mathcal{K}$ is a product of cones, our implementation disaggregates the elements of $\mathcal{K}^*$ per individual cone, adding one OA cut per cone per iteration.

\section{Extended formulations and disciplined convex programming}\label{sec:dcp}

While many problems are representable in conic form, the transformation from the user's algebraic representation of the problem often requires expert knowledge. Disciplined convex programming (DCP) is an algebraic modeling concept proposed by Grant et al.~\cite{DCP}, one of whose original motivations was to provide a means to make these transformations automatic and transparent to users. DCP is not intrinsically tied to conic representations, however. In this section, we present the basic concepts of DCP from the viewpoint of extended formulations. This perspective both provides insight into how conic formulations are generated and enables further generalization of the technique to problems which are not conic representable using standard cones.

Detection of convexity of arbitrary nonlinear expressions is NP-Hard~\cite{ConvexityNPHard}, and since a conic-form representation is a proof of convexity, it is unreasonable to expect a modeling system to be able to reliably generate these representations from arbitrary input. Instead, DCP requires users to construct expressions whose convexity can be proven by simple composition rules. A DCP implementation (e.g., the MATLAB package CVX) provides a library of basic operations like addition, subtraction, norms, square root, square, geometric mean, logarithms, exponential, entropy $x\log(x)$, powers, absolute value, $\max\{x,y\}$, $\min\{x,y\}$, etc. whose curvature (convex, concave, or affine) and monotonicity properties are known. These basic operations are called \textit{atoms}.

All expressions representing the objective function and constraints are built up via compositions of these atoms in such a way that guarantees convexity. For example, the expression $\exp(x^2+y^2)$ is convex and \textit{DCP compliant} because $\exp(\cdot)$ is convex and monotone increasing and $x^2+y^2$ is convex because it is a convex composition (through addition) of two convex atoms. The expression $\sqrt{xy}$ is concave when $x,y\ge 0$ as we noted previously, but not \textit{DCP compliant} because the inner term $xy$ has indefinite curvature. In this case, users must reformulate their expression using a different atom like $geomean(x,y)$. We refer readers to~\cite{DCP,dcpweb} for further introduction to DCP.

An important yet not well-known aspect of DCP is that the composition rules for DCP have a 1-1 correspondence with the existence of extended formulations of epigraphs. For example, suppose $g$ is convex and monotone increasing and $f$ is convex. Then the function $h(x) := g(f(x))$ is convex and recognized as such by DCP. If $E_h := \{ (x,t) : h(x) \le t \}$ is the epigraph of $h$, then we can represent $E_h$ through an extended formulation using the epigraphs $E_g$ and $E_f$ of $g$ and $f$, respectively. That is, $(x,t) \in E_h$ iff $\exists\, s$ such that $(x,s) \in E_f \text{ and } (s,t) \in E_g$.
The validity of this extended formulation follows directly from monotonicty of $g$. Furthermore, if $E_f$ and $E_g$ are conic representable, then so is $E_h$, which is precisely how DCP automatically generates conic formulations. The conic form representation is not necessary, however; one may instead represent $E_f$ and $E_g$ using smooth nonlinear constraints if $f$ and $g$ are smooth.

This correspondence between composition of functions and extended formulations was considered by Tawarmalani and Sahinidis~\cite{Baron}, although in the context of nonconvex MINLP. Composition generalizes the notion of separability far beyond summations of univariate functions as proposed by Hijazi et al.~\cite{Hijazi}. 

DCP, based on the philosophy that users should be ``disciplined'' in their modeling of convex functions, describes a simple set of rules for verifying convexity and rejects any expressions not satisfying them; it is not based on ad-hoc detection of convexity which is common among nonconvex MINLP solvers. DCP is well established within the convex optimization community as a practical modeling technique, and many would agree that it is reasonable to ask users to formulate convex optimization problems in DCP form. By doing so they unknowingly provide all of the information needed to generate powerful extended formulations.

\section{Computational results}

In this section we present preliminary computational results implementing the extended formulations proposed in this work. We have implemented a solver, \textit{Pajarito}, which currently accepts input as mixed-integer conic programming problems with a mix of second-order and exponential cones. We have translated 194 convex problems from MINLPLIB2 representable using these cones into Convex.jl~\cite{convexjl}, a DCP algebraic modeling language in Julia which performs automatic transformation into conic form. We exclude instances without integer constraints, some which are pure quadratic, and some which Bonmin is unable to solve within time limits. Pajarito currently implements traditional OA using derivative-based NLP solvers~\cite{Bonmin} applied to the conic extended formulation, as the conic solvers we tested were not sufficiently reliable. Pajarito itself relies on JuMP~\cite{LubinDunningIJOC}, and the implementation of the core algorithm spans less than 1000 lines of code. Pajarito will be released as open source in the upcoming months.

\begin{figure}[t]
\centering
\subfigure[Solution time]{
\includegraphics[scale=0.28]{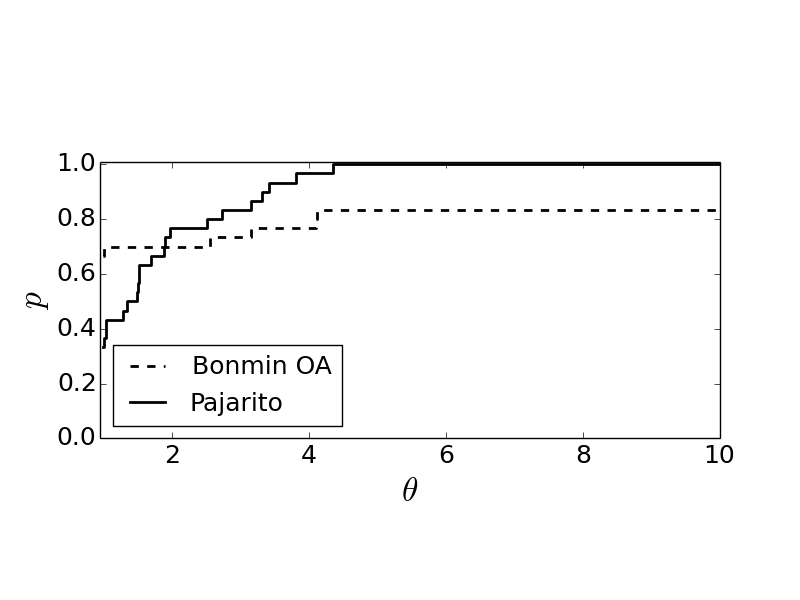}}
\subfigure[Number of OA iterations]{
\includegraphics[scale=0.28]{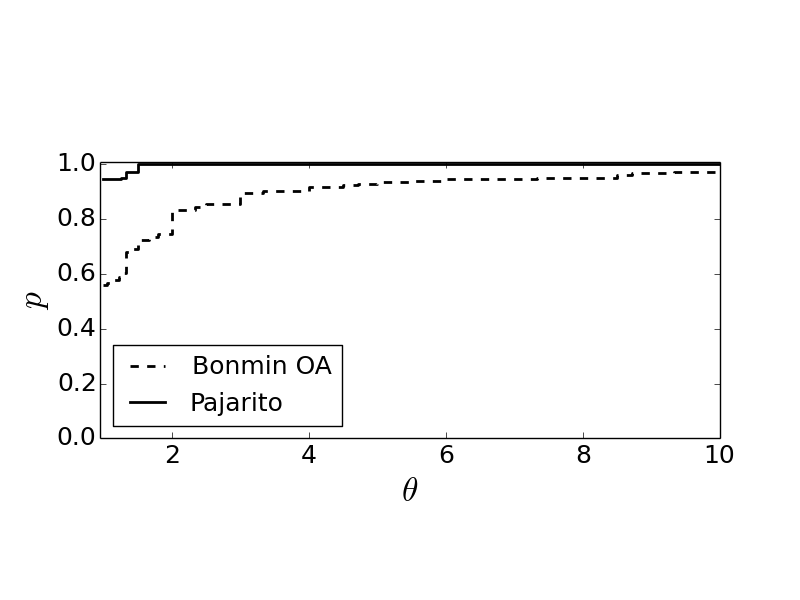}}
\vspace{-.25cm}
\caption{Comparison performance profiles for the entire data set. Higher is better. Here, $p$ is the proportion of instances for which the given solver is within a factor of $\theta$ of the best solution time or iteration count.}
\label{fig:all}
\end{figure}

Our main comparison is with Bonmin's OA algorithm, which in 2014 benchmarks by H. Mittelmann was found to be the overall fastest MICP solver when using CPLEX as the inner MILP solver~\cite{hans}. For the MISOCP instances, we also compare with CPLEX. 
Performance profiles~\cite{perf} of all instances solved by Bonmin in greater than 30 seconds are provided in Figure \ref{fig:all}.  Tables~\ref{tab:results:1} and~\ref{tab:results:2} in the Appendix list the complete set of results. 
Their highlights are:
\begin{enumerate}
    \item We observe that the extended formulation helps significantly reduce the number of OA iterations (Figure~\ref{fig:all}). This can be seen as a sign of scalability provided by the extended formulation.

\item Pajarito is much faster on many of the challenging problems (\texttt{slay},\texttt{netmod},\\ \texttt{portfol\_classical}), although these problems are MISOCPs where CPLEX dominates (note that CPLEX 12.6.2 already applies extended formulations for MISOCPs~\cite{VielmaExtendedFormulations}). Pajarito has not been optimized for performance, leading Bonmin to be faster on the relatively easy instances.
\item Perhaps the strongest demonstration of Pajarito's strength is the \texttt{gams01} instance, which was previously unsolved and whose conic representation requires a mix of SOC and EXP cones. The best known bound was 1735.06 and the best known solution was 21516.83. Pajarito solved the instance to optimality with an objective value of 21380.20 in 6 iterations. Unfortunately, the origin of the instance is unknown and confidential.
\end{enumerate}

\section*{Acknowledgements}
We thank the anonymous referees for their comments. They greatly improved the clarity of the manuscript. We also thank one of the
anonymous referees for pointing out the SOC-representability of the \texttt{sssd} family of instances originally derived in~\cite{perspective}. 
M. Lubin was supported by the DOE Computational Science Graduate
Fellowship, which is provided under grant number DE-FG02-97ER25308.
The work at LANL was funded by the Center for Nonlinear Studies (CNLS) and was carried out under the auspices
of the NNSA of the U.S.
DOE at LANL
under Contract No. DE-AC52-06NA25396. J.P. Vielma was funded by NSF grant CMMI-1351619.

\bibliography{refs}{}

\bibliographystyle{siam}

\newpage

\section*{Appendix}
\vspace{-.25cm}
All computations were performed on a high-performance cluster at Los Alamos with Intel$^\circledR$ Xeon$^\circledR$ E5-2687W v3 @3.10GHz 25.6MB L3 cache processors and 251GB DDR3 memory installed on every node. CPLEX v12.6.2 is used as a MILP and MISOCP solver. We use KNITRO v9.1.0 as an NLP solver for Pajarito. Bonmin v1.8.3 is compiled with CPLEX v12.6.2 and Ipopt 3.12.3 using the HSL linear algebra library MA97. All solvers are set to a relative optimality gap of $10^{-5}$, are run on a single thread (both CPLEX and KNITRO), and are given 10 hours of wall time limit (with the exception of \texttt{gams01} where we give 32 threads to CPLEX for the MILP relaxations).
\vspace{-.7cm}
\begin{figure}[H]
\centering
\subfigure[Solution time]{
\includegraphics[scale=0.28]{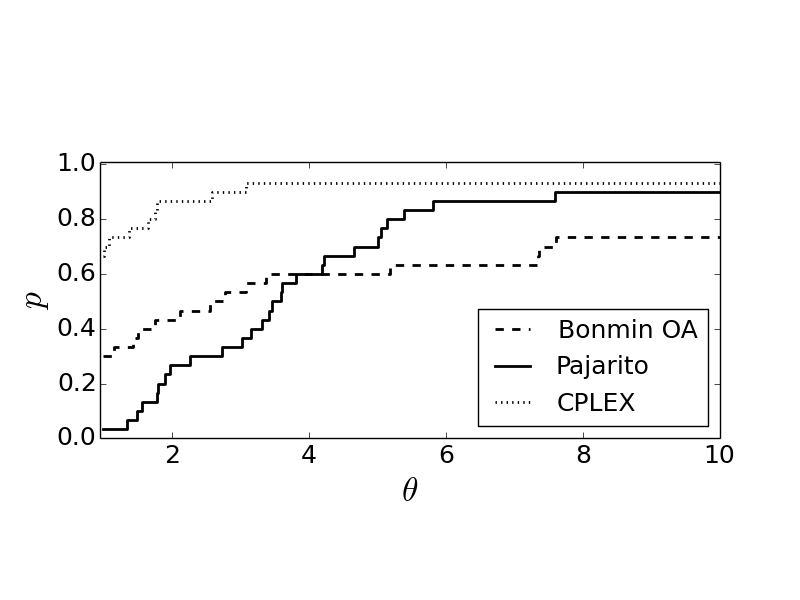}}
\subfigure[Number of OA iterations]{
\includegraphics[scale=0.28]{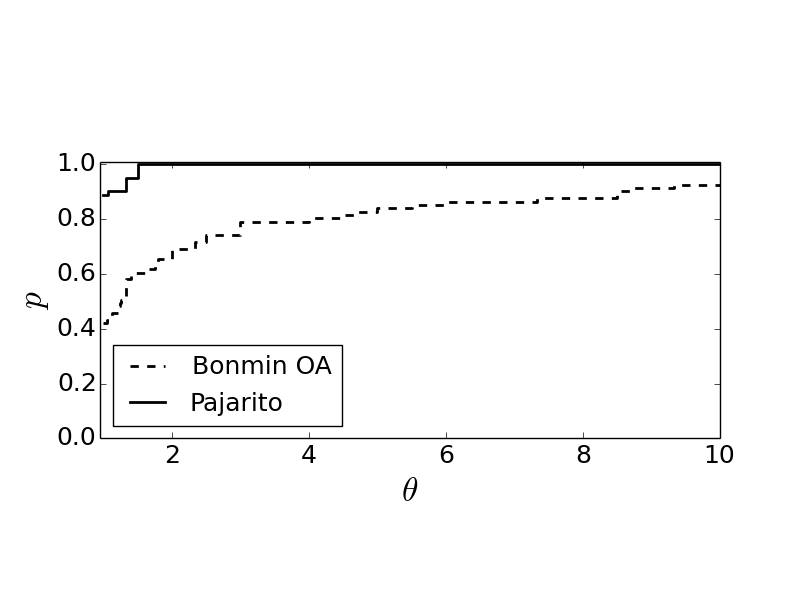}}
\vspace{-.25cm}
\caption{Comparison performance profiles for SOC representable instances}
\label{fig:soc}
\end{figure}
\vspace{-.75cm}
Performance profiles with respect to SOC representable instances are provided in Figure \ref{fig:soc}. After reformulating these instances, we are able solve them using CPLEX as MISOCPs. Although CPLEX dominates, Pajarito is able to solve more instances than Bonmin in this case.

 \begin{table}[t]
 \tiny
 \centering
 \begin{tabular}{l|c|r|r|r|r|r}
   Instance & Conic rep. & Bonmin Iter & Bonmin Time & Pajarito Iter & Pajarito Time & CPLEX Time \\
  \hline
batch & Exp & 2 &     0.60 & 1 &     4.95 & -- \\ 
batchdes & Exp & 1 &     0.07 & 1 &     4.76 & -- \\ 
batchs101006m & Exp & 10 &     1.88 & 3 &     7.67 & -- \\ 
batchs121208m & Exp & 4 &     3.14 & 3 &    11.01 & -- \\ 
batchs151208m & Exp & 6 &     7.97 & 3 &    14.63 & -- \\ 
batchs201210m & Exp & 8 &    14.92 & 2 &    29.09 & -- \\ 
clay0203h & SOC & 9 &     0.90 & 5 &     6.53 &     0.35 \\ 
clay0203m & SOC & 10 &     0.40 & 6 &     6.74 &     0.37 \\ 
clay0204h & SOC & 3 &     3.60 & 1 &     6.27 &     1.61 \\ 
clay0204m & SOC & 3 &     0.33 & 1 &     5.14 &     1.02 \\ 
clay0205h & SOC & 4 &    20.89 & 3 &    28.76 &     8.93 \\ 
clay0205m & SOC & 6 &     5.50 & 3 &    12.67 &     1.77 \\ 
clay0303h & SOC & 9 &     0.97 & 5 &     7.29 &     0.54 \\ 
clay0303m & SOC & 10 &     0.58 & 7 &     7.36 &     0.68 \\ 
clay0304h & SOC & 11 &     5.27 & 9 &    17.96 &     1.42 \\ 
clay0304m & SOC & 16 &     2.84 & 13 &    22.99 &     2.13 \\ 
clay0305h & SOC & 4 &    23.81 & 3 &    56.93 &    23.32 \\ 
clay0305m & SOC & 7 &     6.16 & 3 &    16.14 &     2.51 \\ 
du-opt & SOC & 61 &     0.76 & 7 &     8.69 &     1.54 \\ 
du-opt5 & SOC & 22 &     0.22 & 4 &     6.66 &     1.97 \\ 
enpro48pb & Exp & 2 &     0.22 & 1 &     5.04 & -- \\ 
enpro56pb & Exp & 1 &     0.22 & 1 &     5.11 & -- \\ 
ex1223 & ExpSOC & 3 &     0.07 & 1 &     5.47 & -- \\ 
ex1223a & SOC & 1 &     0.03 & 0 &     4.57 &     0.01 \\ 
ex1223b & ExpSOC & 3 &     0.07 & 1 &     5.52 & -- \\ 
ex4 & SOC & 2 &     0.13 & 2 &     5.80 &     0.86 \\ 
fac3 & SOC & 6 &     0.15 & 2 &     5.22 &     0.07 \\ 
netmod\_dol2 & SOC & 33 &   167.49 & 7 &    53.04 &    12.58 \\ 
netmod\_kar1 & SOC & 102 &    56.45 & 12 &    13.75 &     7.68 \\ 
netmod\_kar2 & SOC & 102 &    56.35 & 12 &    13.68 &     7.66 \\ 
no7\_ar25\_1 & SOC & 2 &    25.19 & 3 &    69.55 &    54.34 \\ 
no7\_ar3\_1 & SOC & 4 &    71.04 & 3 &    95.84 &   126.09 \\ 
no7\_ar4\_1 & SOC & 5 &    85.87 & 4 &   110.70 &    48.97 \\ 
no7\_ar5\_1 & SOC & 7 &    69.23 & 5 &   117.65 &    32.60 \\ 
nvs03 & SOC & 1 &     0.06 & 1 &     4.89 &     0.00 \\ 
slay04h & SOC & 5 &     0.19 & 2 &     5.22 &     0.14 \\ 
slay04m & SOC & 5 &     0.11 & 2 &     5.20 &     0.18 \\ 
slay05h & SOC & 9 &     0.60 & 3 &     5.73 &     0.37 \\ 
slay05m & SOC & 7 &     0.18 & 3 &     5.51 &     0.16 \\ 
slay06h & SOC & 12 &     1.94 & 2 &     5.56 &     0.69 \\ 
slay06m & SOC & 9 &     0.29 & 2 &     5.27 &     0.42 \\ 
slay07h & SOC & 15 &     5.04 & 3 &     6.61 &     0.98 \\ 
slay07m & SOC & 12 &     0.66 & 3 &     5.67 &     0.67 \\ 
slay08h & SOC & 22 &    27.27 & 3 &     7.41 &     1.50 \\ 
slay08m & SOC & 21 &     2.89 & 2 &     5.43 &     0.96 \\ 
slay09h & SOC & 36 &   163.31 & 3 &     9.01 &     1.93 \\ 
slay09m & SOC & 28 &    17.22 & 3 &     6.12 &     1.54 \\ 
slay10h & SOC & 80 &  8155.02 & 4 &    26.13 &     7.55 \\ 
slay10m & SOC & 77 &  1410.08 & 4 &     9.08 &     1.80 \\ 
syn05h & Exp & 2 &     0.09 & 1 &     4.75 & -- \\ 
syn05m & Exp & 2 &     0.07 & 1 &     4.73 & -- \\ 
syn05m02h & Exp & 1 &     0.06 & 1 &     4.79 & -- \\ 
syn05m02m & Exp & 1 &     0.07 & 1 &     4.80 & -- \\ 
syn05m03h & Exp & 1 &     0.07 & 1 &     4.86 & -- \\ 
syn05m03m & Exp & 1 &     0.07 & 1 &     4.83 & -- \\ 
syn05m04h & Exp & 1 &     0.07 & 1 &     4.85 & -- \\ 
syn05m04m & Exp & 1 &     0.08 & 1 &     4.85 & -- \\ 
syn10h & Exp & 1 &     0.04 & 0 &     4.46 & -- \\ 
syn10m & Exp & 2 &     0.04 & 1 &     4.79 & -- \\ 
syn10m02h & Exp & 1 &     0.09 & 1 &     4.92 & -- \\ 
syn10m02m & Exp & 2 &     0.09 & 1 &     4.85 & -- \\ 
syn10m03h & Exp & 1 &     0.08 & 1 &     4.91 & -- \\ 
syn10m03m & Exp & 1 &     0.08 & 1 &     4.85 & -- \\ 
syn10m04h & Exp & 1 &     0.11 & 1 &     5.03 & -- \\ 
syn10m04m & Exp & 1 &     0.11 & 1 &     5.01 & -- \\ 
syn15h & Exp & 1 &     0.06 & 1 &     4.86 & -- \\ 
syn15m & Exp & 2 &     0.07 & 1 &     4.79 & -- \\ 
syn15m02h & Exp & 1 &     0.09 & 1 &     5.00 & -- \\ 
syn15m02m & Exp & 1 &     0.09 & 1 &     4.92 & -- \\ 
syn15m03h & Exp & 1 &     0.13 & 1 &    48.61 & -- \\ 
syn15m03m & Exp & 2 &     0.11 & 1 &     4.94 & -- \\ 
syn15m04h & Exp & 1 &     0.14 & 1 &     5.77 & -- \\ 
syn15m04m & Exp & 2 &     0.14 & 1 &     5.10 & -- \\ 
syn20h & Exp & 2 &     0.10 & 2 &     5.14 & -- \\ 
syn20m & Exp & 2 &     0.06 & 1 &     4.81 & -- \\ 
syn20m02h & Exp & 2 &     0.15 & 2 &     5.70 & -- \\ 
syn20m02m & Exp & 2 &     0.10 & 2 &     5.24 & -- \\ 
syn20m03h & Exp & 1 &     0.13 & 1 &     5.55 & -- \\ 
syn20m03m & Exp & 2 &     0.15 & 2 &     5.34 & -- \\ 
syn20m04h & Exp & 1 &     0.19 & 1 &     6.03 & -- \\ 
syn20m04m & Exp & 2 &     0.27 & 2 &     5.60 & -- \\ 
syn30h & Exp & 3 &     0.12 & 3 &     5.61 & -- \\ 
syn30m & Exp & 3 &     0.09 & 3 &     5.40 & -- \\ 
syn30m02h & Exp & 3 &     0.21 & 3 &     6.34 & -- \\ 
syn30m02m & Exp & 4 &     0.19 & 3 &     5.69 & -- \\ 
syn30m03h & Exp & 3 &     0.40 & 3 &     6.86 & -- \\ 
syn30m03m & Exp & 3 &     0.27 & 3 &     6.16 & -- \\ 
syn30m04h & Exp & 3 &     0.49 & 3 &     7.91 & -- \\ 
syn30m04m & Exp & 4 &     0.42 & 3 &     6.60 & -- \\ 
syn40h & Exp & 4 &     0.19 & 3 &     5.74 & -- \\ 
syn40m & Exp & 4 &     0.97 & 2 &     5.20 & -- \\ 
syn40m02h & Exp & 3 &     0.31 & 3 &     6.74 & -- \\ 
syn40m02m & Exp & 3 &     0.24 & 3 &     5.97 & -- \\ 
syn40m03h & Exp & 4 &     0.59 & 4 &     8.66 & -- \\ 
syn40m03m & Exp & 5 &     0.52 & 4 &     7.17 & -- \\ 
syn40m04h & Exp & 4 &     1.02 & 4 &    10.42 & -- \\ 
syn40m04m & Exp & 5 &     0.87 & 5 &     9.25 & -- \\ 
 \end{tabular}
 \vspace{0.3cm}
 \caption{MINLPLIB2 instances. ``Conic rep'' column indicates which cones are used in the conic representation of the instance (second-order cone and/or exponential). CPLEX is capable of solving only second-order cone instances. Times in seconds.}
 \label{tab:results:1}
 \end{table}

 \begin{table}[t]
 \tiny
 \centering
 \begin{tabular}{l|c|r|r|r|r|r}
   Instance & Conic rep. & Bonmin Iter & Bonmin Time & Pajarito Iter & Pajarito Time & CPLEX Time \\
  \hline
synthes1 & Exp & 3 &     0.04 & 2 &     5.04 & -- \\ 
synthes2 & Exp & 3 &     0.05 & 2 &     4.98 & -- \\ 
synthes3 & Exp & 6 &     0.10 & 2 &     5.00 & -- \\ 
rsyn0805h & Exp & 1 &     0.14 & 1 &     4.92 & -- \\ 
rsyn0805m & Exp & 2 &     0.25 & 2 &     5.22 & -- \\ 
rsyn0805m02h & Exp & 5 &     0.71 & 5 &     7.31 & -- \\ 
rsyn0805m02m & Exp & 4 &     2.16 & 4 &     7.24 & -- \\ 
rsyn0805m03m & Exp & 3 &     4.08 & 3 &     7.76 & -- \\ 
rsyn0805m04m & Exp & 2 &     2.31 & 2 &     6.78 & -- \\ 
rsyn0810m & Exp & 2 &     0.24 & 1 &     4.94 & -- \\ 
rsyn0810m02h & Exp & 3 &     0.58 & 3 &     6.45 & -- \\ 
rsyn0810m02m & Exp & 4 &     5.78 & 3 &     6.83 & -- \\ 
rsyn0810m03h & Exp & 3 &     1.36 & 3 &     7.62 & -- \\ 
rsyn0810m03m & Exp & 3 &     6.04 & 3 &     8.66 & -- \\ 
rsyn0810m04h & Exp & 3 &     1.31 & 2 &     7.71 & -- \\ 
rsyn0810m04m & Exp & 4 &     3.77 & 3 &     8.14 & -- \\ 
rsyn0815h & Exp & 1 &     0.27 & 1 &    23.50 & -- \\ 
rsyn0815m & Exp & 2 &     0.23 & 2 &     5.25 & -- \\ 
rsyn0815m02m & Exp & 5 &     1.94 & 4 &     7.14 & -- \\ 
rsyn0815m03h & Exp & 5 &     5.21 & 5 &    16.04 & -- \\ 
rsyn0815m03m & Exp & 4 &     4.59 & 5 &    10.16 & -- \\ 
rsyn0815m04h & Exp & 3 &     2.03 & 3 &    10.43 & -- \\ 
rsyn0815m04m & Exp & 4 &     7.78 & 3 &    10.68 & -- \\ 
rsyn0820h & Exp & 3 &     0.42 & 2 &     5.59 & -- \\ 
rsyn0820m & Exp & 2 &     0.24 & 2 &     5.29 & -- \\ 
rsyn0820m02h & Exp & 3 &     0.59 & 2 &     6.72 & -- \\ 
rsyn0820m02m & Exp & 3 &     1.90 & 3 &     6.90 & -- \\ 
rsyn0820m03h & Exp & 2 &     1.37 & 2 &     7.76 & -- \\ 
rsyn0820m03m & Exp & 3 &     5.14 & 3 &     8.83 & -- \\ 
rsyn0820m04h & Exp & 4 &     2.66 & 4 &    11.59 & -- \\ 
rsyn0820m04m & Exp & 3 &     8.65 & 3 &    11.52 & -- \\ 
rsyn0830h & Exp & 3 &     0.41 & 3 &     5.95 & -- \\ 
rsyn0830m & Exp & 4 &     0.37 & 4 &     6.19 & -- \\ 
rsyn0830m02m & Exp & 5 &     1.83 & 5 &    15.68 & -- \\ 
rsyn0830m03h & Exp & 2 &     1.45 & 2 &     9.04 & -- \\ 
rsyn0830m03m & Exp & 4 &     3.45 & 4 &    10.15 & -- \\ 
rsyn0830m04h & Exp & 3 &     2.35 & 3 &    12.59 & -- \\ 
rsyn0830m04m & Exp & 4 &    11.47 & 4 &    15.82 & -- \\ 
rsyn0840h & Exp & 2 &     0.30 & 2 &     5.69 & -- \\ 
rsyn0840m & Exp & 2 &     0.26 & 3 &     5.72 & -- \\ 
rsyn0840m02h & Exp & 3 &     0.72 & 2 &     7.34 & -- \\ 
rsyn0840m02m & Exp & 4 &     1.53 & 3 &     7.73 & -- \\ 
rsyn0840m03h & Exp & 3 &     1.85 & 3 &    11.07 & -- \\ 
rsyn0840m03m & Exp & 5 &     2.47 & 5 &    12.41 & -- \\ 
rsyn0840m04h & Exp & 2 &     2.40 & 2 &    44.19 & -- \\ 
rsyn0840m04m & Exp & 4 &     7.62 & 4 &    22.33 & -- \\ 
sambal & SOC & 0 &     0.03 & 0 &     4.52 &     0.00 \\ 
gbd & SOC & 1 &     0.04 & 0 &     4.55 &     0.00 \\ 
ravempb & Exp & 4 &     0.33 & 1 &     5.22 & -- \\ 
portfol\_classical050\_1 & SOC & \textgreater 989 & \textgreater 36000 & 12 &    37.77 &     3.31 \\ 
m3 & SOC & 1 &     0.68 & 0 &     4.58 &     0.07 \\ 
m6 & SOC & 1 &     0.16 & 1 &     5.18 &     0.17 \\ 
m7 & SOC & 1 &     0.59 & 0 &     4.84 &     0.69 \\ 
m7\_ar25\_1 & SOC & 1 &     0.37 & 1 &     5.19 &     0.16 \\ 
m7\_ar2\_1 & SOC & 1 &     2.19 & 1 &     7.01 &     1.58 \\ 
m7\_ar3\_1 & SOC & 1 &     1.88 & 1 &     6.79 &     0.82 \\ 
m7\_ar4\_1 & SOC & 1 &     0.35 & 0 &     4.77 &     0.84 \\ 
m7\_ar5\_1 & SOC & 1 &     0.34 & 0 &     5.71 &     0.98 \\ 
fo7 & SOC & 3 &    27.68 & 4 &    42.88 &    23.67 \\ 
fo7\_2 & SOC & 2 &    12.52 & 2 &    16.70 &     4.88 \\ 
fo7\_ar25\_1 & SOC & 4 &     9.87 & 4 &    27.18 &     9.92 \\ 
fo7\_ar2\_1 & SOC & 2 &     8.68 & 3 &    19.63 &    11.04 \\ 
fo7\_ar3\_1 & SOC & 3 &    11.61 & 3 &    31.28 &    22.16 \\ 
fo7\_ar4\_1 & SOC & 2 &     9.61 & 2 &    15.68 &    10.27 \\ 
fo7\_ar5\_1 & SOC & 1 &     5.66 & 1 &     8.95 &    12.67 \\ 
fo8 & SOC & 2 &    79.50 & 3 &    82.41 &    52.92 \\ 
fo8\_ar25\_1 & SOC & 3 &    45.80 & 4 &   144.43 &    63.09 \\ 
fo8\_ar2\_1 & SOC & 3 &    59.24 & 4 &   161.68 &    60.09 \\ 
fo8\_ar3\_1 & SOC & 1 &    14.65 & 1 &    14.78 &    37.85 \\ 
fo8\_ar4\_1 & SOC & 1 &    10.53 & 1 &    16.48 &    62.60 \\ 
fo8\_ar5\_1 & SOC & 2 &    23.26 & 1 &    34.09 &    59.75 \\ 
fo9 & SOC & 3 &   534.56 & 4 &   209.68 &   227.52 \\ 
fo9\_ar25\_1 & SOC & 6 &  1430.17 & 6 &  6221.39 &  1240.89 \\ 
fo9\_ar3\_1 & SOC & 1 &    16.77 & 1 &    22.69 &   103.84 \\ 
fo9\_ar4\_1 & SOC & 2 &    40.77 & 1 &    60.73 &   785.75 \\ 
fo9\_ar5\_1 & SOC & 2 &    39.47 & 3 &   134.95 &   725.60 \\ 
flay02h & SOC & 2 &     0.09 & 2 &     5.18 &     0.02 \\ 
flay02m & SOC & 2 &     0.05 & 2 &     5.12 &     0.04 \\ 
flay03h & SOC & 8 &     0.40 & 8 &     7.08 &     0.20 \\ 
flay03m & SOC & 8 &     0.17 & 8 &     6.74 &     0.24 \\ 
flay04h & SOC & 24 &    19.92 & 24 &    30.22 &     1.14 \\ 
flay04m & SOC & 22 &     4.43 & 22 &    16.03 &     1.00 \\ 
flay05h & SOC & 181 &  6583.08 & 164 &  6593.05 &    96.62 \\ 
flay05m & SOC & 180 &  3258.45 & 171 &  4938.36 &    68.91 \\ 
flay06h & SOC & \textgreater 30 & \textgreater 36000 & \textgreater 32 & \textgreater 36000 &  6958.36 \\ 
flay06m & SOC & \textgreater 68 & \textgreater 36000 & \textgreater 55 & \textgreater 36000 &  4752.04 \\ 
o7 & SOC & 9 &  1623.33 & 8 &  3060.63 &   526.94 \\ 
o7\_2 & SOC & 5 &   435.47 & 5 &   663.47 &   128.95 \\ 
o7\_ar25\_1 & SOC & 4 &   259.10 & 3 &   510.12 &   455.29 \\ 
o7\_ar2\_1 & SOC & 1 &    41.51 & 1 &   137.82 &    68.66 \\ 
o7\_ar3\_1 & SOC & 4 &   338.68 & 3 &   642.90 &   875.63 \\ 
o7\_ar4\_1 & SOC & 7 &  1486.87 & 7 &  2239.11 &   535.17 \\ 
o7\_ar5\_1 & SOC & 4 &   309.86 & 4 &   777.35 &   216.84 \\ 
o8\_ar4\_1 & SOC & 4 &  2736.05 & 3 & 10438.68 &  8447.35 \\ 
tls2 & SOC & 7 &     0.19 & 4 &     5.27 &     0.10 \\ 
tls4 & SOC & 88 &   260.67 & 7 &    18.58 &     6.15 \\ 
{\bf gams01} & {\bf ExpSOC} & {\bf \textgreater 19} & {\bf \textgreater 36000} & {\bf 6} & {\bf 23414.37} & -- \\
 \end{tabular}
 \vspace{0.3cm}
 \caption{MINLPLIB2 instances, continued.}
 \label{tab:results:2}
 \end{table}

\end{document}